\documentclass[12pt,leqno]{amsart}

\usepackage[utf8]{inputenc}   
\usepackage[T1]{fontenc}      
\usepackage{geometry}         
\usepackage[all]{xy}
\usepackage{amssymb}
\usepackage{enumerate}
\usepackage{enumitem}
\usepackage{stmaryrd}

\voffset-1cm \markleft{\rm} 
\newtheorem{theorem}{Theorem}
\newtheorem{definition}[theorem]{Definition}
\newtheorem{corollary}[theorem]{Corollary}

\newtheorem{lemma}[theorem]{Lemma}

\newtheorem{proposition}[theorem]{Proposition}

\newtheorem{defi}{D\'{e}finition}

\newcommand{\K}{\mathbb K}
\newcommand{\C}{\mathbb C}

\newcommand{\R}{\mathbb R}

\newcommand{\ds}{\displaystyle}

\usepackage{graphicx}

\title{Decomposition  of trilinear forms}
\author{Elisabeth Remm. Michel Goze}
\date{}
\address{Elisabeth Remm Universit\'e de Haute-Alsace, IRIMAS UR 7499, F-68100 Mulhouse, France.}
\email{elisabeth.remm@uha.fr}
\email{goze.rac@gmail.com}

\begin{document}

\maketitle

\begin{abstract}
In this work, we are interested in the decomposition of  trilinear forms. We associate this decomposition with the decomposition into irreducible subspaces of the group algebra associated with the symmetric group. We have a correspondence between a ternary form and a symmetric trilinear form, in the same way that we have the correspondence between a quadratic form and a symmetric quadratic form. We are then interested in the decomposition of ternary forms into a sum of cubes of linear forms and as an application we find the Tartaglia formulae on the expression of the roots of equations of degree 3.
\end{abstract}




\bigskip

\noindent{\bf Introduction.}
The decomposition of a bilinear form into a sum of a symmetric form and a skew-symmetric form is actually based on the decomposition into irreducible subspaces of the associative group algebra defined by the symmetric group $\Sigma_2$ of degree $2$. Starting from this remark, we propose a decomposition of trilinear forms by describing a decomposition of the group algebra  associated with the symmetric group $\Sigma_3$ of degree $3$. Recall that if $G$ is a finite group, the group algebra  $A(G)$ defined by $G$ is the finite dimensional vector space generated by $G$ and this space is naturally endowed with an associative algebra structure and a $G$-module structure. 

We extend this result to trilinear applications with values in a vector space. This is the case, for example, of the associator of a multiplication of a non-associative algebra. This makes it possible to highlight particular classes of such algebras.

To complete this work, we are interested in ternary forms (also called cubic forms) associated with symmetric trilinear forms. We show under which conditions these ternary forms on a two-dimensional $\K$-vector space $E$ can be decomposed into a sum of two cubes of linear forms. As a direct consequence, we find again the Tartaglia formulae giving the expressions of roots of real equations of degree $3$.

We briefly touch upon finding solutions to degree 5 equations using this method. 

The trilinear forms and more generally the multilinear forms occur in different mathematical areas and subject such as cryptography (see \cite{B.S.}) and also in the study  of properties of algebra laws (see \cite{B.B.R.}).

\section{Definition}
\subsection{Trilinear forms}

Let $E$ be a $\K$-vector space where  $\K=\R$ or $\C$. A map
$$f: E^3 \rightarrow \K$$
is a trilinear form if $f$  is a linear form with respect to each of its arguments. If $\dim E=n$ and if $\{e_1,\cdots,e_n\}$ is a basis of $E$, then the writing of $f$ in this basis is an homogeneous algebraic expression of degree $3$. The constants
$$a_{ijk}=f(e_i,e_j,e_k)$$
are the structure constants (or coefficients) of $f$ with respect to this basis. The formulae of change of basis are the classical formulae corresponding to $f \in (E \otimes E \otimes E)^ \ast$.
\begin{definition}
Let $f$ be a trilinear form on $E$. Then
\begin{enumerate}
\item $f$ is symmetric if for every $x_1,x_2,x_3 \in E$,
$$\forall \sigma \in \Sigma_3, \ \ f(x_{\sigma(1)},x_{\sigma(2)},x_{\sigma(3)})=f(x_1,x_2,x_3),$$
\item $f$ is alternating if it vanishes when two inputs are the same.
\end{enumerate}
\end{definition}
A trilinear form on $E$ is antisymmetric (skew-symmetric) if
$$ f(x_{\sigma(1)},x_{\sigma(2)},x_{\sigma(3)})=-f(x_1,x_2,x_3)$$
for every transposition $\sigma \in\Sigma_3$ and for every $x_1,x_2,x_3 \in E.$
But, since we assume that the characteristic of $\K$ is $0$, the notions of antisymmetric and alternating trilinear forms coincide.

A $p$-linear form $f:E^p \rightarrow \K$ is defined in the same way as a linear form with respect to each of its $p$ arguments.

\subsection{Symmetrization-alternatization }
Let $f$ be a trilinear form on $E$. The trilinear form $f_s$ defined by
$$f_s(x_1,x_2,x_3)=\sum_{\sigma \in \Sigma_3}f(x_{\sigma(1)},x_{\sigma(2)},x_{\sigma(3)})$$
is symmetric. It is called the \textbf{symmetrization} of $f$.

\medskip

The trilinear form $f_a$ defined by
$$f_a(x_1,x_2,x_3)=\sum_{\sigma \in \Sigma_3}\mathrm{sgn(\sigma)}f(x_{\sigma(1)},x_{\sigma(2)},x_{\sigma(3)})$$
is alternating. It is called the \textbf{alternatization} of $f$.

\medskip\textbf{Remark : case $p=2$.} It is well known, that for bilinear forms we have a converse to this construction. Given a symmetric bilinear form $f_1$ and an alternating bilinear form $f_2$, there is a unique bilinear form $f$ such that 
$$f_s=f_1, \ \ f_a=f_2.$$
This form is given by
$$f=\frac{f_s+f_a}{2}$$
(recall that we assume that $\emph{Char}(\K) =0).$ 

In the next section, we will give an interpretation  of this decomposition in terms of  representations of $\Sigma_3$ to extend it to the trilinear case.
\section{The group algebra $\K[\Sigma_p]$}
\subsection{Generalities}
If $\Sigma_p$ is the symmetric group of degree $p$, then the group algebra $\K[\Sigma_p]$  is the associative algebra over $\K$ 
whose underlying $\K$-module is the the free module over $\K$ on the underlying set of $\Sigma_p$ and 
whose multiplication is given on basis elements by the group operation. An element of $\K[\Sigma_p]$ is a formal combination
$$v=\sum_{\sigma\in \Sigma_p}a_\sigma \sigma.$$
The  action of $\Sigma_p$ on $\K[\Sigma_p]$ is defined by
$$(\sigma, v=\sum_{i}a_i \sigma_i) \longrightarrow \sigma v=\sum_{i}a_i \sigma^{-1} \circ \sigma_i.$$ 
For each $v \in \K[\Sigma_p]$, we denote by $\mathcal{O}(v)$ the corresponding orbit and by $F_v =
\K(\mathcal{O}(v))$ the linear subspace of $\K[\Sigma_p]$ generated by $\mathcal{O}(v)$. A subspace $F$ is called \textbf{invariant} if for every $\sigma \in \Sigma_p$,  we have $\sigma F \subset F$. In particular the subspaces $F_v$ are invariant.

\medskip

\noindent\textbf{Examples. }
We denote by $\tau_{i,j}$ the transposition between $i$ and $j$. 
\begin{enumerate}
\item $p=2$. We have
$$\Sigma_2=\{Id,\tau_{1,2}\}$$
and
$$\K[\Sigma_2]=\{a_1 Id+a_2\tau_{1,2}, \ a_i \in \K \}.$$
If $v=a_1Id+a_2\tau_{1,2}$, then $\mathcal{O}(v)=\{v, a_2 Id+a_1\tau_{1,2}\}.$ We deduce that $\dim F_v=2$ except if $a_1^2-a_2^2=0.$
\begin{proposition}
Assume that $v \neq 0.$ Then $\dim F_v=2$ except if $v=a_1(Id-\tau_{1,2})$ or $v=a_1(Id +\tau_{1,2})$. In these last cases, if $a_1 \neq 0$, $\dim F_v=1$
\end{proposition}
If we denote by
$$V_2=Id-\tau_{1,2}, \ \ W_2=Id+\tau_{1,2}$$
then  $F_{V_2}$ and $F_{W_2}$ are the only one-dimensional invariant subspaces of $\K[\Sigma_2]$ and we have the classical decomposition
$$\K[\Sigma_2]=F_{V_2} \oplus F_{W_2}.$$
\item $p=3$. We have
$$\Sigma_3=\{Id, \tau_{1,2}, \tau_{1,3} , \tau_{2,3}, c, c^2\}$$
where $c$ is the cycle $(2,3,1)$ and $c^2=c \circ c.$ If we denote by $V_3$ and $W_3$ the vectors
$$V_3=Id-\tau_{1,2}-\tau_{1,3} - \tau_{2,3}+ c+ c^2, \ \ W_3=Id+\tau_{1,2}+\tau_{1,3} + \tau_{2,3}+ c+ c^2,$$
then $\dim F_{V_3}=\dim F_{W_3}=1$ and these spaces are  the only one-dimensional invariant subspaces of $\K[\Sigma_3]$.
\end{enumerate}
\subsection{Irreducible subspaces  of dimension $2$ for $p=3$}
An invariant subspace $F$ of $\K[\Sigma_3]$ is called \textbf{irreducible} if any invariant subspace of $F$ is trivial. The only one-dimensional irreducible subspaces are $F_{V_3}$ and $F_{W_3}$. 

Assume now that $F$ is a $2$-dimensional irreducible subspace of $\K[\Sigma_3]$. For any $v \neq 0 \in F$,  $\{\sigma v, \sigma \in \Sigma_3\}$ is a generating family of $F$ and there exists $\sigma \in \Sigma_3$, such as $\{v,\sigma v\}$ is a basis of $F$.

Let $v=a_1Id +a_2 \tau _{12} +a_3  \tau _{13} +a_4 \tau _{23}+a_5 c_1 +a_6c_2$ be
in $\mathbb{K}\,[\Sigma_3].$ 
\subsubsection{The family $\{v,\tau_{1,2}v\}$ is linearly dependent}
The hypothesis implies $\tau_{1,2}v=\lambda v$. Solving this linear system, we obtain $\lambda=\pm 1$ and
  $v=(a_1,a_1,a_3,a_4,a_4,a_3)$ or $(a_1,-a_1,a_3,a_4,-a_4,-a_3)$.  Since $\dim F_v=2$, then there exists $\sigma \in \Sigma_3$ such that $\{v, \sigma v\}$ is a basis of $F_v$. If we also have that $\tau_{1,3}v=\lambda v$, then $v=V_3$ or $W_3$ and $F_V$ is not irreducible. Then we can assume that $\{v,\tau_{1,3}v\}$ is a basis of $F_v$.

  In the first case, that is $v=(a_1,a_1,a_3,a_4,a_4,a_3)$, we have
  $$v+\tau_{1,3}v+\tau_{2,3}v=(a_1+a_3+a_4)W_3.$$
  If $a_1+a_3+a_4 \neq 0$, $W_3\in F_v$ and this subspace is not irreducible. Then $a_1+a_3+a_4=0$ and 
  $$v=(a_1,a_1,a_3,-a_1-a_3,-a_1-a_3,a_3).$$
  
  In the second case, that is $v=(a_1,-a_1,a_3,a_4,-a_4,-a_3)$, we have
  $$v-\tau_{1,3}v-\tau_{2,3}v=(a_1-a_3-a_4)V_3$$
  This implies $a_4=a_1-a_3$ and
  $$v=(a_1,-a_1,a_3,a_1-a_3,-a_1+a_3,-a_3).$$
  
\subsubsection{The family $\{v,\tau_{1,2}v\}$ is a basis of $F_v$}  
  
In this case we consider the vector $u=v+\tau_{1,2}v$. This vector is in $F_v$ then $F_u$ is a subspace of $F_v$. If $\dim F_u=1$, then $F_u=F_{V_3}$ or $F_{W_3}$ and $F_v$ is not irreducible. Then $F_u=F_v$. But $\{u,\tau_{1,2}u\}$ is linearly dependent.
So we return to the previous case.

\begin{proposition}
Let us consider the vectors
$$u_1=(1,1,0,-1,-1,0), \ \ u_2=(1,1,-1,0,0,-1).$$ 
and
$$u_3=(1,-1,0,1,-1,0), \ \ u_4=(1,-1,1,0,0,-1).$$ 
Then the subspaces $F_{u_i}$  are $2$-dimensional irreducible subspaces. 
\end{proposition}

\begin{theorem} 

We have the following direct decomposition of $\K[\Sigma_3]$ in irreducible subspaces
$$\K[\Sigma_3]=F_{V_3} \oplus F_{W_3} \oplus F_{u_i} \oplus F_{u_j}, \ i\neq j.$$
\end{theorem}
In fact the vectors $V_3,W_3,u_i,\tau_{13}u_i,u_j,\tau_{13}u_j, \ (i\neq j)$ are linearly independent.

\section{Decomposition of a trilinear form}
Let $\varphi$ be a $p$-linear form on $E$ and $v=\sum a_i\sigma_i \in \K[\Sigma_p]$. We define the $p$-linear form $\varphi_v$ by
$$\varphi_v(x_1,\cdots,x_p)=\sum a_i \varphi(x_{\sigma_i(1)},\cdots,x_{\sigma_i(p)}).$$
For example, the symmetric and alternating part of $\varphi$ correspond to
$$\varphi_s=\varphi_{V_p}, \ \varphi_a=\varphi_{W_p}.$$
In particular, if $\varphi$ is a bilinear form, the decomposition in symmetric and alternating part corresponds to
$$\varphi=\frac{1}{2}(\varphi_{V_2}+\varphi_{W_2})$$
and it is associated with the decomposition
$$\R[\Sigma_2]=F_{V_2} \oplus F_{W_2}$$
in irreducible subspaces.
We generalize this remark to the trilinear case. 
\subsection{Decomposition of trilinear maps}
We consider the following basis of $\R[\Sigma_3]$ 
$$\mathcal{B}=\{V_3,W_3,u_1,u_2,u_3,u_4\}.$$
\begin{lemma}
We have 
$$6 Id=V_3+W_3+u_1+u_2+u_3+u_4.$$
\end{lemma}
\begin{proof}
In fact
$$
\begin{array}{l}
V_3=Id-\tau_{1,2}-\tau_{1,3} - \tau_{2,3}+ c+ c^2=(1,-1,-1,-1,1,1), \\ 
W_3=Id+\tau_{1,2}+\tau_{1,3} + \tau_{2,3}+ c+ c^2=(1,1,1,1,1,1), \\
u_1=Id+\tau_{1,2}-\tau_{2,3}-c=(1,1,0,-1,-1,0), \\
u_2=Id+\tau_{1,2}-\tau_{1,3}-c^2=(1,1,-1,0,0,-1), \\
u_3=Id-\tau_{1,2}+\tau_{2,3}-c=(1,-1,0,1,-1,0),\\
u_4=Id-\tau_{1,2}+\tau_{1,3}-c^2=(1,-1,1,0,0,-1)\\
\end{array}
$$
We deduce the identity.
\end{proof}
We define the polarisation of a trilinear form
\begin{definition}
Let $\varphi$ be a real trilinear form. The polarization of $\varphi$ is the decomposition
$$\varphi=\frac{1}{6}(\varphi_s+\varphi_s+\varphi_{u_1}+\varphi_{u_1}+\varphi_{u_2}+\varphi_{u_3}+\varphi_{u_4})$$
\end{definition}

So we have
$$
\begin{array}{l}
\varphi_a(x,y,z) =  \varphi(x,y,z)-\varphi(y,x,z)-\varphi(z,y,x)-\varphi(x,z,y)+\varphi(y,z,x)+\varphi(z,x,y),\\
\varphi_s(x,y,z) =  \varphi(x,y,z)+\varphi(y,x,z)+\varphi(z,y,x)+\varphi(x,z,y)+\varphi(y,z,x)+\varphi(z,x,y),\\
\varphi_{u_1}(x,y,z)= \varphi(x,y,z)+\varphi(y,x,z)-\varphi(x,z,y)-\varphi(y,z,x),\\
\varphi_{u_2}(x,y,z)= \varphi(x,y,z)+\varphi(y,x,z)-\varphi(z,y,x)-\varphi(z,x,y).\\
\varphi_{u_3}(x,y,z)= \varphi(x,y,z)-\varphi(y,x,z)+\varphi(x,z,y)-\varphi(y,z,x),\\
\varphi_{u_4}(x,y,z)= \varphi(x,y,z)-\varphi(y,x,z)+\varphi(z,y,x)-\varphi(z,x,y).\\
\end{array}
$$
\begin{proposition}
The trilinear forms $\varphi_{u_1}$ and  $\varphi_{u_2}$ satisfy the following properties: 

$\varphi_{u_i}(x,y,z)-\varphi_{u_i}(y,x,z)=0, \ \
\varphi_{u_i}(x,y,z)+\varphi_{u_i}(z,y,x)+\varphi_{u_i}(x,z,y)=0.$

The trilinear forms $\varphi_{u_3}$ and  $\varphi_{u_4}$ satisfy the following properties: 

$\varphi_{u_i}(x,y,z)+\varphi_{u_i}(y,x,z)=0, \ \
\varphi_{u_i}(x,y,z)-\varphi_{u_i}(z,y,x)-\varphi_{u_i}(x,z,y)=0.$

\end{proposition}

\subsection{Polynomial expression}
We assume $E=\K^n$. Let $\varphi$ be a trilinear form on $E$. If
$$x_1=(x_1^1,x_1^2,\cdots,x_1^n), x_2=(x_2^1,x_2^2,\cdots,,x_2^n), x_3=(x_3^1,x_3^2,\cdots,x_3^n)$$
are three vectors of $E$, then
$$\varphi(x_1,x_2,x_3)=\sum a_{i,j,k}x_1^{i}x_2^{j}x_3^{k}$$
is identified with an homogeneous polynomial expression $P_\varphi$ of degree 3 in the variables $x_i^j$.

1. If $\varphi=\varphi_s$ is symmetric, then $P\varphi$ is a symmetric polynomial. But any symmetric polynomial in $n^3$ variables
$x_1^1,x_1^2,\cdots,x_1^n,x_2^1,x_2^2,\cdots,,x_2^n,x_3^1,x_3^2,\cdots,x_3^n$
 is representable in a unique way as a polynomial in the elementary symmetric polynomials. But $P_\varphi$ is homogeneous of degree $3$. The previous decomposition is quite trivial  but we can refine this decomposition. 
 
 \noindent\textbf{Examples}
 
 i) $n=2$ We have the decomposition
 $$\varphi(x_1,x_2,x_3)=a_1x_1^1x_2^1x_3^1 +a_2x_1^2x_2^2x_3^2+a_3(x_1^1x_2^1x_3^2+x_1^1x_2^2x_3^1+x_1^2x_2^1x_3^1)+a_4(x_1^1x_2^2x_3^2+x_1^2x_2^2x_3^1+x_1^2x_2^1x_3^2).$$
 
ii)  $n=3$, we have 

\begin{proposition}
Any symmetric trilinear form on $E$ is a linear combination of the following symmetric forms:
\begin{enumerate}
\item $\varphi_{s,1,2,3}(x_1,x_2,x_3)=\sum_{i\neq j \neq k} x_1^i x_2^jx_3^k$ \\
\item $\varphi_{s,i,j}(x_1,x_2,x_3)=x_1^ix_2^ix_3^j+x_2^ix_3^ix_1^j+x_1^ix_3^ix_2^j$, \ $i \neq j$ \\
\item $\varphi_{s,i}(x_1,x_2,x_3)=x_1^ix_2^ix_3^i,\ i=1,2,3.$
\end{enumerate}
\end{proposition}

2. If $\varphi$ is a trilinear alternating form, we have the decomposition
$$\varphi(x_1,x_2,x_3)= \sum_{1 \leq i <j<k\leq n}a_{ijk}\sum_{\sigma \in \Sigma_3} \epsilon(\sigma) x_1^{\sigma(i)}x_2^{\sigma(j)}x_3^{\sigma(k)}
$$
For example, if $n=3$ then $\varphi$ is colinear to the trilinear form corresponding to the determinant, that is
$$\varphi_a(x_1,x_2,x_3)=\det (x_1,x_2,x_3)=\sum_{\sigma \in \Sigma_3}\varepsilon(\sigma)x_1^{\sigma(1)}x_2^{\sigma(2)}x_3^{\sigma(3)}.$$

\medskip

\noindent From now, we assume that $n=3$. In this case we have 
$\varphi(x_1,x_2,x_3)=P_\varphi$ is a polynomial in variables $x_1^1,x_2^1,x_3^1,x_1^2,x_2^2,x_3^2,x_1^3,x_2^3,x_3^3.$

3. We saw that the form $\varphi_{u_1}$ satisfies the following symmetry relations:
$$
\begin{array}{l}
\varphi_{u_1}(x_1,x_2,x_3)-\varphi_{u_1}(x_2,x_1,x_3)=0,\\
\varphi_{u_1}(x_1,x_2,x_3)+\varphi_{u_1}(x_3,x_2,x_1)+\varphi_{u_1}(x_1,x_3,x_2)=0.\\
\end{array}
$$
We have
\begin{proposition}
The polynomial associated with 
$\varphi_{u_1}$ or $\varphi_{u_2}$ is a linear combination of the following polynomials
$$
\left\{
\begin{array}{l}
 2x_1^ix_2^ix_3^j-x_1^ix_2^jx_3^i-x_1^jx_2^ix_3^i  \ \ i\neq j,\\
x_1^1x_2^2x_3^3+x_1^2x_2^1x_3^3-x_1^1x_2^3x_3^2-x_1^3x_2^1x_3^2,\\
x_1^2x_3^3x_3^1+x_1^3x_2^2x_3^1-x_1^1x_2^3x_3^2-x_1^3x_2^1x_3^2.\\
\end{array}
\right.
$$
\end{proposition}

4. The form $\varphi_{u_3}$ (and $\varphi_{u_4}$) satisfies the following relations:
$$
\begin{array}{l}
\varphi_{u_3}(x_1,x_2,x_3)+\varphi_{u_3}(x_2,x_1,x_3)=0, \\
\varphi_{u_3}(x_1,x_2,x_3)-\varphi_{u_3}(x_3,x_2,x_1)-\varphi_{u_3}(x_1,x_3,x_2)=0.
\end{array}
$$
We have
\begin{proposition}
The polynomial associated with $\varphi_{u_3}$ or $\varphi_{u_4}$ is a linear combination of the following polynomials
$$
\left\{
\begin{array}{l}
 x_1^ix_2^jx_3^i-x_1^jx_2^ix_3^i  \ \ i\neq j,\\
x_1^1x_2^2x_3^3-x_1^2x_2^1x_3^3+x_1^1x_2^3x_3^2-x_1^3x_2^1x_3^2,\\
x_1^2x_2^3x_3^1-x_1^3x_2^2x_3^1+x_1^1x_2^3x_3^2-x_1^3x_2^1x_3^2.\\
\end{array}
\right.
$$
\end{proposition}
\section{Nonassociative algebras. A decomposition of the associator.}
The multiplication  $\mu$ of an algebra $A$ over $\K$ is given by a bilinear map
$$\mu: A \times A \rightarrow A.$$
Its associator $\mathcal{A}_\mu$ is the trilinear map given by $\mathcal{A}_\mu(x,y,z)=\mu(\mu(x,y,z)-\mu(x,\mu(y,z)).$ For simplification, we will denote $xy $ instead of $\mu(x,y)$ and refer to the associator as $\mathcal{A}$. The algebra $(A,\mu)$ is associative when $\mathcal{A}=0$. It is nonassociative if $\mathcal{A}$ is null or not null. As with trilinear forms, for any vector $v=\sum a_i\sigma_i \in \K[\Sigma_3]$, we define the trilinear map $\mathcal{A}_v$ by
$$\mathcal{A}_v(x_1,x_2,x_3)=\sum a_i\mathcal{A}(x_{\sigma_i(1)},x_{\sigma_i(2)},x_{\sigma_i(3)}).$$
According to the previous section, we have
$$\mathcal{A}=\frac{1}{6}(\mathcal{A}_{V_3}+\mathcal{A}_{W_3}+\mathcal{A}_{u_1}+\mathcal{A}_{u_2}+\mathcal{A}_{u_3}+\mathcal{A}_{u_4}).$$
This decomposition highlights the following  classes of nonassociative algebras:

1. \textit{Totally symmetric algebras}. The associator $\mathcal{A}_{V_3}$ satisfies the relations
$$\mathcal{A}_{V_3}(x,y,z)=\mathcal{A}_{V_3}(y,x,z)=\mathcal{A}_{V_3}(z,y,x).$$ If $(A,\mu)$ is a nonassociative algebra whose associator corresponds to $\mathcal{A}_{V_3}$, it satisfies the previous symmetric relations and it is a nonassocative totally symmetric algebra.

Remark: If we consider a nonassociave algebra whose associator $\mathcal{A}$ satisfies
$$\mathcal{A}_{V_3}=0,$$
that is,
$$\mathcal{A}(x_1,x_2,x_3)+\mathcal{A}(x_2,x_1,x_3)+\mathcal{A}(x_3,x_2,x_1)+\mathcal{A}(x_1,x_3,x_2)+\mathcal{A}(x_2,x_3,x_1)+\mathcal{A}(x_3,x_1,x_2)=0$$
it is a  $3$-power associative algebras.

2. \textit{Alternative algebras}. The associator $\mathcal{A}_{W_3}$ satisfies the relations
$$\mathcal{A}_{W_3}(x,y,z)=-\mathcal{A}_{W_3}(y,x,z)=-\mathcal{A}_{W_3}(z,y,x).$$ If $(A,\mu)$ is a nonassociative algebra whose associator corresponds to $\mathcal{A}_{W_3}$, it satisfies the previous symmetric relations and it is an alternative algebra.

\noindent Remark: If we consider a nonassociave algebra whose associator $\mathcal{A}$ satisfies
$$\mathcal{A}_{W_3}=0$$
that is
$$\mathcal{A}(x_1,x_2,x_3)-\mathcal{A}(x_2,x_1,x_3)-\mathcal{A}(x_3,x_2,x_1)-\mathcal{A}(x_1,x_3,x_2)+\mathcal{A}(x_2,x_3,x_1)+\mathcal{A}(x_3,x_1,x_2)=0$$
it is a Lie-admissible algebra.

Numerous studies have focused on these two classes. See for example \cite{G.R.}.

3. \textit{ Left-symmetric $3$-power associative algebras.}
Since $$\mathcal{A}_{u_1}(x_1,x_2,x_3)=\mathcal{A}(x_1,x_2,x_3)+\mathcal{A}(x_2,x_1,x_3)-\mathcal{A}(x_1,x_3,x_2)-\mathcal{A}(x_2,x_3,x_1)$$
this trilinear map satisfies
$$
\left\{
\begin{array}{l}
\mathcal{A}_{u_1}(x_1,x_2,x_3)-\mathcal{A}_{u_1}(x_2,x_1,x_3)=0,\\
\mathcal{A}_{u_1}(x_1,x_2,x_3)+\mathcal{A}_{u_1}(x_3,x_2,x_1)+\mathcal{A}_{u_1}(x_1,x_3,x_2)=0
\end{array}
\right.
$$
An algebra $A$ whose associator satisfies the first relation  is left symmetric which implies that it is also Lie-admissible  and it is $3$-power associative  from the second relation. 

Remark. Let us consider the class of nonassociative  algebras whose associator $\mathcal{A}$ satisfies
$$\mathcal{A}_{u_1}=0$$
that is
$$\mathcal{A}(x_1,x_2,x_3)+\mathcal{A}(x_2,x_1,x_3)-\mathcal{A}(x_1,x_3,x_2)-\mathcal{A}(x_2,x_3,x_1)=0.$$
This class does not seem to have been the subject of particular studies. However, let us note that this class contains the class of right-symmetric  algebras.

4.  \textit{Left-alternative Lie-admissible algebra.}
We have 
$$\mathcal{A}_{u_3}(x_1,x_2,x_3)=\mathcal{A}(x_1,x_2,x_3)-\mathcal{A}(x_2,x_1,x_3)+\mathcal{A}(x_1,x_3,x_2)-\mathcal{A}(x_2,x_3,x_1)$$
these trilinear map satisfies
$$
\left\{
\begin{array}{l}
\mathcal{A}_{u_3}(x_1,x_2,x_3)+\mathcal{A}_{u_3}(x_2,x_1,x_3)=0,\\
\mathcal{A}_{u_3}(x_1,x_2,x_3)-\mathcal{A}_{u_3}(x_3,x_2,x_1)-\mathcal{A}_{u_3}(x_1,x_3,x_2)=0
\end{array}
\right.
$$
The first relation shows that the algebra 
is left-alternative which implies that it is $3$-power-associative. The second identity implies that this algebra is Lie-admissible.

\noindent Remark. If we consider a nonassociave algebra whose associator $\mathcal{A}$ satisfies
$$\mathcal{A}_{u_3}=0,$$
that is,
$$\mathcal{A}(x_1,x_2,x_3)-\mathcal{A}(x_2,x_1,x_3)+\mathcal{A}(x_1,x_3,x_2)-\mathcal{A}(x_2,x_3,x_1)=0.$$
This class of nonassociative algebras contains the class of right-alternative algebras that is satisfying
$$\mathcal{A}(x_1,x_2,x_3)+\mathcal{A}(x_1,x_3,x_2).$$

\section{On the reduction of cubic forms}

\subsection{Symmetric trilinear forms and cubic forms}
The decomposition of a trilinear form focuses on  the study of  antisymmetric, symmetric forms and those defined by vectors $u_1$ and $u_2$.
Let $E$ be a finite-dimensional $\K$-vector space. The classification of alternating trilinear forms is already presented in \cite{Cap}. The study of bilinear symmetric forms is associated with the study of quadratic forms. 

Recall that a quadratic form  is a map $q : E \rightarrow \K$  such that $q(av) = a^2q(v)$ for all $a \in \K, \ v \in  E$ and the function $f:E\times E \rightarrow \K$ defined by $f(u,v)=\frac{1}{2}(q(u+v)-q(u) - q(v))$ is a bilinear form on $E$. Such a map can be identified to a homogeneous polynomial of degree $2$ in $n$ variables with coefficients in $\K$.
A major result specifies that any quadratic form can be written as a linear combination of squares of linearly independent linear forms. What about the cubic forms? We assume now that $\K=\C$. The problem that naturally arises is that of the decomposition of cubic forms.

   A cubic form is a homogeneous polynomial of degree $3$ in $n$ variables: For example
\begin{itemize}
\item \emph{binary} cubic form: $q(x,y)=ax^3+bx^2y+cxy^2+dy^3$
\item \emph{ternary} cubic form: $q(x,y,z)=ax^3+bx^2y+cx^2z+dxy^2+exz^2+fxyz+gy^3+hy^2z+kyz^2+lz^3$
\end{itemize}
with $a,b,\cdots,l \in \K$.

If $\varphi$ is a trilinear  form on the $\K$ vector space $E$, the associated cubic form is 
$$q_{\varphi}(u)=\varphi (u,u,u).$$
It is clear that  if $\varphi$ is skewsymmetric then $ q_{\varphi}=0$. But from Proposition 7,  it is the same if $\varphi=\varphi_{u_i}$ for $i \in \{ 1,2,3,4\}$. We deduce that the cubic forms are in a one-to-one correspondence with the symmetric trilinear forms. 

\noindent The next section concerns the reduction  of binary cubic forms, that is when $E$ is 
$2$-dimensional. 

\subsection{Decomposition of a binary cubic form}
In this paragraph, we are interested in the decomposition of a  binary cubic form into a sum of cubes of linear forms. 
\begin{proposition}
There are binary cubic forms that cannot be written as a sum of cubes of independent linear forms.
\end{proposition}
\begin{proof}
Let us consider, for example, the form
$$q(x,y)=x^3+x^2y.$$
Assume that
$$q(x,y)=(\alpha_1x+\beta_1y)^3+(\alpha_2x+\beta_2y)^3.$$
This is equivalent to the algebraic system
$$
\left\{
\begin{array}{l}
\alpha_1^3+\alpha_2^3=1\\
3\alpha_1^2 \beta_1+3\alpha_2^2\beta_2=1 \\
3\alpha_1\beta_1^2+3\alpha_2\beta_2^2=0 \\
\beta_1^3+\beta_2^3=0
\end{array}
\right.
$$
Let us note that $\beta_2=0$ implies $\beta_1=0$ and $\alpha_1^2\beta_1+\alpha_2^2\beta_2=0$ which is impossible. Then $\beta_2 \neq 0$ (and also $\alpha_1,\alpha_2,\beta_1$). We deduce
$$\alpha_2=-\frac{\alpha_1\beta_1^2}{\beta_2^2}$$
and
$$3(\alpha_1^2 \beta_1+\alpha_2^2\beta_2)=3\left(\alpha_1^2 \beta_1+\frac{\alpha_1^2\beta_1^4}{\beta_2^4}\beta_2\right)=3\alpha_1^2\beta_1\left(1+\frac{\beta_1^3}{\beta_2^3}\right)=3\alpha_1^2\beta_1\left(\frac{\beta_2^3+\beta_1^3}{\beta_2^3}\right)=0\neq 1$$
so this algebraic system has no solution.

\end{proof}
Since there are binary cubic 
forms that are the sum of cubes of linear forms, we can wonder about the structure of such forms that are written in this way.
Let $q(x,y)=ax^3+bx^2y+cxy^2+dy^3$ a binary cubic form. When $a\neq 0$, a classical change of variables $x=X-\frac{b}{3a}y$  leads to $aX^3+pXy^2+qy^3$ with $p=c-\frac{b^2}{3a}, q=d-\frac{bc}{3a}+\frac{2b^3}{27a^2}$ so allows to reduce the problem for $a\neq 0$ to write as the sum of cubes of linear forms the binary cubic forms
$$q(x,y)=ax^3+cxy^2+dy^3 \ { \rm with} \ a\neq 0$$
Assume that
$$q(x,y)=(\alpha_1x+\beta_1y)^3+(\alpha_2x+\beta_2 y)^3.$$
Then we have
$$
\left\{
\begin{array}{l}
\alpha_1^3+\alpha_2^3=a\\
3\alpha_1^2 \beta_1+3\alpha_2^2\beta_2=0 \\
3\alpha_1\beta_1^2+3\alpha_2\beta_2^2=c \\
\beta_1^3+\beta_2^3=d
\end{array}
\right.
$$
If $c=0$ then $q$ is trivially a sum of linear forms.
If one of the coefficients $\alpha_i,\beta_i$ is zero, then $c=0$ and $q$ is indeed a sum of cubes of independent linear forms (two cubes if $d\neq 0$).
Then we assume that $\alpha_1\alpha_2\beta_1\beta_2 \neq 0.$ 
Then $$\alpha_1=\frac{c}{3a}\frac{\alpha_2^3}{\beta_1^2} \ \ \beta_2=-\frac{\alpha_1^2}{\alpha_2^2}\beta_1=-\frac{c^2}{9a^2}\frac{\alpha_2^4}{\beta_1^3} $$
We put $u=\alpha_2, \ v= \frac{\alpha_2}{\beta_1}.$ Then
$$\alpha_1=\frac{c}{3a}uv^2,\  \beta_2=-\frac{c^2}{9a^2}uv^3, \ \beta_1=\frac{u}{v}. $$
The equations $\alpha_1^3+\alpha_2^3=a, \ \beta_1^3+\beta_2^3=d$ are equivalent to
$$\ds{u^3\left(1+ \frac{c^3}{3^3a^3} v^6\right)=a}, \ \ $$
$$u^3\left(1-\frac{c^6}{9^3a^6}v^{12}\right)=dv^3.$$
We deduce 
$$1-\frac{c^3}{27a^3}v^6=\frac{d}{a}v^3.$$
We put $V=v^3$. Then $V$ is a root of
$$\frac{c^3}{27a^3}V^2+\frac{d}{a}V-1=0.$$
The discriminant is $\Delta=\frac{27ad^2+4c^3}{27a^3}.$ We deduce
$$V=v^3=\frac{27a^3}{2c^3}\left(-\frac{d}{a}\pm \sqrt{\Delta}\right)$$
and $$u^3\left(2-\frac{d}{a}V \right)=a.$$
Since we have assumed that $a \neq 0$, the system has a solution if and only if  
$$ 4c^3+27ad^2 \neq 0$$
and this condition is equivalent to $\Delta\neq 0.$

We can now calculate the coefficients $\alpha_i$ and $\beta_i$. 
We have
$$v^3= \frac{27a^3}{2c^3}\left(-\frac{d}{a}+\epsilon \sqrt{\Delta}\right), \ \ u^3=\frac{a^2}{2a-dv^3} $$
We deduce
$$u^3= \frac{2c^3}{27a\left(a\Delta-\epsilon d\sqrt{\Delta}\right)}$$
with $\epsilon^2=1,$
and
 
$$u=\frac{c}{3}\left(\frac{2}{a^2(a\Delta-\epsilon d \sqrt{\Delta})}\right)^{\frac{1}{3}}, 
\ \ v=\frac{3a}{c}\left(\frac{\epsilon a\sqrt{\Delta}- d}{2a}\right)^{\frac{1}{3}}. $$
We check that the equations  $\alpha_1^3+\alpha_2^3=a$ and $\beta_1^3+\beta_2^3=d$ are indeed satisfied.  Let’s look at the second equation. The first one is identical. We have
$$\alpha_2=u, \ \beta_1=\frac{u}{v}, \ \alpha_1=\frac{c}{3a}uv^2, \ \beta_2=-\frac{c^2}{9a^2}uv^3.$$
We have
$$\beta_1^3+\beta_2^3=\left(\frac{u}{v}\right)^3\left[1-\frac{c^6}{9^3a^6}v^{12}\right].$$
We have seen that
$$\frac{c¨^3}{27a^3}v^6+\frac{d}{a}v^3-1=0.$$
We deduce
$$v^6=\frac{27a^3}{c^3}-\frac{27a^2d}{c^3}v^3.$$
This implies
$$v^{12}=\frac{27^2a^6}{ac^6}\left[ac^3+27a^2d^2-\left(27ad^3+2dc^3\right) v^3\right].$$
We put
$$\rho=ac^3+27a^2d^2-(27ad^3+2dc^3)v^3$$
and $\beta_1^3+\beta_2^3=d$ is equivalent to
$$u^3\left(1-\frac{\rho}{ac^3}\right)=dv^3.$$
But $u^3=\frac{a^2}{2a-dv^3}$. Then we have to verify that
$$a^2\left(1-\frac{\rho}{ac^3}\right)=2adv^3-d^2v^6=2adv^3-d^2\left(\frac{27a^3}{c^3}-\frac{27a^2d}{c^3}v^3\right).$$
This is equivalent to
$$ac^3-\rho=(2dc^3+27ad^3)v^3-27a^2d^2.$$
And we find again the expression of $\rho$ and the equation $\beta_1^3+\beta_2^3=d$ is satisfied.

\medskip

\begin{theorem}
Let $q(x,y)=ax^3+cxy^2+dy^3$ be a ternary form with $a\neq 0$. 

If $c=0,$ the binary cubic form $q(x,y)=ax^3+dy^3$ writes trivially as a sum of two linear forms.

If $c\neq 0$   we can write $q$ as a sum of the cube of two linear forms
if and only if $$\delta= 4c^3+27ad^2 \neq 0$$ and we have
$$q(x,y)=(\alpha_1x+\beta_1y)^3+(\alpha_2x+\beta_2y)^3$$
with
$$\alpha_2=u, \ \beta_1=\frac{u}{v}, \ \alpha_1=\frac{c}{3a}uv^2, \ \beta_2=-\frac{c^2}{9a^2}uv^3$$
and $$u=\frac{c}{3}\left(\frac{2}{a^2(a\Delta-d\sqrt{\Delta}} \right) ^\frac{1}{3}, \ \ v=\frac{3a}{c} \left( -\frac{d}{2a}+\frac{\sqrt{\Delta}}{2} \right) ^\frac{1}{3}$$
with $\Delta= \frac{\delta}{27a^3}.$
\end{theorem}

It  remains to examine the case 
 $q(x,y)=bx^2y+cxy^2+dy^3$, that is $a=0$ in the initial case . 
 If $d\neq 0$ the study is similar to the previous case exchanging the role of $x$ and $y$.
  We can therefore only study the particular case
$$q(x,y)=bx^2y+cxy^2.$$
If $bc=0$, then $q$ cannot be written as a sum of the cube of linear forms. Assume $bc\neq 0$. Since $\alpha_1^3+\alpha_2^3=\beta_1^3+\beta_2^3=0$, we can choose $\alpha_1=k_1\alpha_2, \beta_1=k_2\beta_2$ with $k_1^3=k_2^3=-1$ and $1+k_1^2k_2 \neq 0$ and $1+k_1k_2^2 \neq 0.$ Let's take, for example, $k_1=-1,k_2=-j$. In this case $1+k_1^2k_2=1-j$ and $1+k_1k_2^2=2+j.$ We deduce
$$\alpha_2^3=\frac{b^2(1+k_1k_2^2)}{3c(1+k_1^2k_2)^2}=\frac{b^2(1+2j)}{9c}, \ \beta_2^3=-\frac{(1+2j)c^2}{3b}.$$
\begin{proposition}
If $bc \neq 0$ then
$$\begin{array}{l}
q(x,y) =  bx^2y+cxy^2 \\
\qquad =   \ds{\left[-\left(\frac{b^2(1+2j)}{9c}\right)^{\frac{1}{3}}\! \! x+j^2\! \! \left(\frac{(1+2j)c^2}{3b}\right)^\frac{1}{3}\! \! y\right]^3\! \! +\! \! \left[\left(\frac{b^2(1+2j)}{9c}\right)^{\frac{1}{3}}\! \! x-j \! \! \left(\frac{(1+2j)c^2}{3b}\right)^\frac{1}{3}\! \! y\right]^3}.
\end{array}$$
\end{proposition}

\noindent{\textbf{Example}.  We saw that $q(x,y)=x^3+x^2y$ could not decompose like this. The change of coordinates $x_1=x-\frac{1}{3}y,y=y$  allows  to write in canonical form
$$q(x_1,y)=x_1^3-\frac{1}{3}x_1y ^2+\frac{2}{27}y^3$$
and $\Delta=0.$ But in this case, although $q$ does not decompose as a sum of cubes, we have the decomposition 
$$q(x_1,y)=\left(x_1-\frac{1}{3}y\right)^3+\left(x_1-\frac{1}{3}y\right)^2y.$$
We will see that when $\delta=27ad^2+4c^3=0,$ $q$ will decompose in this form, that is $q(x,y)=(\alpha_1x+\beta_1y)^3+(\alpha_1x+\beta_1y)^2(\alpha_2x+\beta_2y).$

\begin{proposition}
Let be $q(x,y)=ax^3+cxy^2+dy^3$ with $a \neq 0$ and $\delta=27ad^2+4c^3=0.$ Then we have the following decomposition
$$q(x,y)=(\alpha_1x+\beta_1y)^3+(\alpha_1x+\beta_1y)^2(\alpha_2x+\beta_2y).$$
\end{proposition}
\begin{proof}
The algebraic system associated with this decomposition is
$$
\left\{
\begin{array}{l}
\alpha_1^3+\alpha_1^2\alpha_2=a\\
3\alpha_1^2\beta_1+\alpha_1^2\beta_2+2\alpha_1\alpha_2\beta_1=0\\
3\alpha_1\beta_1^2+2\alpha_1\beta_1\beta_2+\alpha_2\beta_1^2=c\\
\beta_1^3+\beta_1^2\beta_2=d.
\end{array}
\right.
$$
Let us note that if $c=0$, $d=0$ and $q$ is already decomposed. Then $c\neq 0.$ In this case the coefficients $\alpha_1,\beta_1$ are non zero. We deduce
$$
\alpha_2=\frac{a-\alpha_1^3}{\alpha_1^2},\ \ 
\beta_1=\frac{3d}{2c}\alpha_1,\ \ 
\beta_2=\frac{4c^2}{9d\alpha_1^2}-\frac{3d}{2c}\alpha_1.
$$
\end{proof}
In particular, if $\alpha_1=-a^\frac{1}{3}$, then  $\alpha_2=2a^\frac{1}{3}$ and $\beta_1=-\frac{3d}{2c}a^\frac{1}{3}=\beta_2.$
\begin{corollary}
Consider $q(x,y)=ax^3+cxy^2+dy^3$ with $a \neq 0$ and $27ad^2+4c^3=0.$ Then we have the following decomposition
$$q(x,y)=\left(-a^\frac{1}{3}x-\frac{3d}{2c}a^\frac{1}{3}y\right)^3+\left(-a^\frac{1}{3}x-\frac{3d}{2c}a^\frac{1}{3}y\right)^2\left(2a^\frac{1}{3}x-\frac{3d}{2c}a^\frac{1}{3}y\right).$$
\end{corollary}
This is equivalent to
$$q(x,y)=-a\left[\left(x-\frac{3d}{2c}y\right)^3-\left(x-\frac{3d}{2c}y\right)^2\left(2x-\frac{3d}{2c}a^\frac{1}{3}y\right)\right].$$

\medskip

\section{Roots of  the equation $ax^3+bx^2+cx+d=0.$}
To find the roots of an equation of degree 3, we have only to consider the case $a \neq 0$.
Using the previous decompositions, we can also consider $b=0$, that is, $ax^3+cx+d=0.$
If $c(4c^3+27ad^2) = 0$, the roots are easy to find.
If $c(4c^3+27ad^2) \neq 0$, we have
$$q(x,1)=ax^3+cx+d=(\alpha_1x+\beta_1)^3+(\alpha_2x+\beta_2)^3$$
with
$$\alpha_2=u, \ \beta_1=\frac{u}{v}, \ \alpha_1=\frac{c}{3a}uv^2, \ \beta_2=-\frac{c^2}{9a^2}uv^3$$
and $$u=\frac{c}{3}\left(\frac{2}{a^2(a\Delta-d\sqrt{\Delta}}\right)^\frac{1}{3}, \ \ v=\frac{3a}{c}\left(-\frac{d}{2a}+\frac{\sqrt{\Delta}}{2}\right)^\frac{1}{3}$$
We deduce that
$$x_1= -\frac{\beta_1+\beta_2}{\alpha_1+\alpha_2}$$
is a root of the equation $ax^3+cx+d=0.$
\begin{theorem}(\textit{Tartaglia's roots}) \ 
If $c(4c^3+27ad^2) \neq 0$ and $a \neq 0$ then
$$x_1=-\frac{1}{v}+\frac{c}{3a}v$$
with $$
v=\frac{3a}{c}\left(\frac{-d+\sqrt{\Delta}}{2a}\right)^\frac{1}{3}$$
is a root of $ax^3+cx+d=0.$
\end{theorem}
\begin{proof}
$$q(x,y)=(\alpha_1x+\beta_1)^3+(\alpha_2x+\beta_2)^3$$
with $\alpha_2=u, \ \beta_1= \frac{u}{v}$, $\alpha_1=\frac{c}{3a}uv^2,\  \beta_2=-\frac{c^2}{9a^2}uv^3$. 
We have a root 
$$x_1=-\frac{\beta_1+\beta_2}{\alpha_1+\alpha_2}=-\frac{1}{v}+\frac{c}{3a}v.$$
\end{proof}

\medskip

\noindent \textbf{Remarks.} If $\Delta > 0$, then $v \in \R$. We have seen that
$$v^3=\frac{27a^3}{2c^3}\left(-\frac{d}{a}\pm \sqrt{\Delta}\right)$$
To find a root of $ax^3+cx+d=0$ we have considered that $\pm =+$.  If $v_1$ corresponds to $:\pm=+$ and $v_2$ to $\pm=-$, then we have $\frac{c}{3a}v_2=-\frac{1}{v_1}$ and we find the same root.

If $\Delta <0$ and $x_1\in \R$, then we find again the previous conditions. If $x_1 \in \C-\R$, then we find the conjugate root.

\noindent

\section{On the equations of degree 5}

In the same vein, we conclude this work by examining what this degree 3 method can induce for degree 5 equations.

Let us consider a
binary pentatic form
$$f(x,y)=ax^5+bx^4y+cx^3y^2+dx^2y^3+exy^4+fy^5.$$
As for the binary cubic case, if $a\neq 0$, we can always consider that $b=0$.  We moreover  assume that $e=0$, that is 
$$f(x,y)=ax^5+cx^3y^2+dx^2y^3+fy^5.$$

If  $f$ is the sum of two $5$-power linear equations:
$$f(x,y)=(\alpha_1x+\beta_1y)^5+(\alpha_2x+\beta_2y)^5,$$
 we have the following algebraic system
$$
\left\{
\begin{array}{l}
\alpha_1^5+\alpha_2^5=a,\\
\alpha_1^4\beta_1+\alpha_2^4\beta_2=0,\\
\alpha_1^3\beta_1^2+\alpha_2^3\beta_2^2=\frac{c}{10},\\
\alpha_1^2\beta_1^3+\alpha_2^2\beta_2^3=\frac{d}{10},\\
\alpha_1\beta_1^4+\alpha_2\beta_2^4=0,\\
\beta_1^5+\beta_2^5=f.
\end{array}
\right..
$$
If $(c,d)=(0,0)$ it is trivial that $f$ is the sum of two $5$-power linear equations.
If $(c,d)\neq (0,0),$ we have $\alpha_1\alpha_2\beta_1\beta_2 \neq 0.$ The second equation is equivalent to
$$\beta_2=-\beta_1\left( \frac{\alpha_1}{\alpha_2}\right)^4.$$
This implies
$$
\left\{
\begin{array}{l}
\alpha_1^3\beta_1^2+\alpha_2^3\beta_2^2=\alpha_1^3\beta_1^2\frac{a}{{\alpha_2}^5}=\frac{c}{10}, \\
\alpha_1^2\beta_1^3+\alpha_2^2\beta_2^3= \alpha_1^2\beta_1^3 \frac{a(\alpha_2^5-\alpha_1^5)}{\alpha_2^{10}}=\frac{d}{10},\\
\alpha_1^{15}+\alpha_2^{15}=0,\\
\alpha_1^5+\alpha_2^5=a,\\
\beta_1^5+\beta_2^5=f.
\end{array}
\right.
$$

But $\alpha_1^{15}+\alpha_2^{15}=(\alpha_1^5+\alpha_2^5)(\alpha_1^{10}-\alpha_1^5\alpha_2^5+\alpha_2^{10})$
implying  $\alpha_1^{10}-\alpha_1^5\alpha_2^5+\alpha_2^{10}=0$ because $a \neq 0$. In this case 
$$\left(\frac{\alpha_1}{\alpha_2}\right)^5=-j \ {\mathrm{or}} \ -j^2.$$
When $\left(\frac{\alpha_1}{\alpha_2}\right)^5=-j^2$ we have $$\alpha_1^5=\frac{-j^2a}{1-j^2}, \ \ \alpha_2^5=\frac{a}{1-j^2}$$
But $\left(\frac{\alpha_1}{\alpha_2}\right)^5=-j$ implies $$\alpha_2^5=\frac{a}{1-j}=\frac{-j^2a}{1-j^2}, \ \ \alpha_1^5=\frac{a}{1-j^2}$$
and all the equations   are symmetric in $\alpha_1,\alpha_2$ so we can consider $\left(\frac{\alpha_1}{\alpha_2}\right)^5=-j^2$ 
 and $$\alpha_1^5=\frac{-j^2a}{1-j^2}, \ \ \alpha_2^5=\frac{a}{1-j^2},\ \ \beta_1^5=\frac{f}{1-j^2}, \ \  \beta_2^5=\frac{-j^2f}{1-j^2}.$$
We have
$$ \frac{\alpha_1}{\alpha_2}=-j\xi,$$
where $\xi$ is a fifth root of the unity.
This implies
$$\alpha_2^3\beta_1^2=\frac{c}{10\xi^3(j^2-1)}, \ \ \alpha_2^2\beta_1^3=\frac{d}{10\xi^2(j^2-1)}.$$
We deduce
$$10^2af=cd$$ and 
$$c^5=-10^5a^3f^2.$$

\begin{theorem}
Let $f(x,y)=ax^5+cx^3y^2+dx^2y^3+fy^5$ a penlinear form with
$$10^2af=cd, c^5=-10^5a^3f^2.$$

Then this form can be written as follows
$$f(x,y)=(\alpha_1 x+\beta_1 y)^5+(\alpha_2 x+\beta_2 y)^5$$
with 
$$\alpha_1=\left(\frac{-j^2a}{1-j^2}\right)^{\frac{1}{5}}, \ \ \alpha_2=\left(\frac{a}{1-j^2}\right)^{\frac{1}{5}}, \ \ \beta_1=\left(\frac{f}{1-j^2}\right)^{\frac{1}{5}}, \ \  \beta_2=\left(\frac{-j^2f}{1-j^2}\right)^{\frac{1}{5}}.$$
\end{theorem}
\begin{corollary}
The equation of degree $5$
$$ax^5+cx^3+dx^2+f=0$$
with 
$$10^2af=cd, c^5=-10^5a^3f^2.$$
admits the following root
$$x=\frac{-\beta_1-\beta_2}{\alpha_1+\alpha_2}.$$
\end{corollary}

\bigskip

\textit{This work was presented at the University of Cagliari in honor of the retirement of Professor Paola Piu.}


\begin{thebibliography}{99}

\bibitem{B.B.R.} Barreiro E., Benayadi S., Rizzo C. {\it Nearly associative algebras}. Journal of Algebra {\bf 658} (2024) 821--868.

\bibitem{B.S.} Boneh D., Silverberg A. {\it Applications of multilinear forms to cryptography}. Contemporary Mathematics {\bf 324}, AMS (2003) 71--90.


\bibitem{Cap} Capdevielle B.   {\it Classification des formes trilinéaires alternées en dimension 6}. L'Enseignement Mathématique {\bf 18}, 1972-1973


\bibitem{G.R.} Goze M., Remm E. 
{\it A Class of Nonassociative Algebras}.
Algebra Colloquium {\bf 14} : 2 (2007) 313–326.

\end{thebibliography}
\end{document}